\documentclass[11pt]{amsart}
\usepackage{amssymb}
\usepackage{latexsym}
\theoremstyle{plain}
\newtheorem{theorem}{Theorem}
\newtheorem{corollary}{Corollary}

 \newtheorem{lemma}{Lemma}
\newtheorem{proposition}{Proposition}

\theoremstyle{definition}

\theoremstyle{remark}

\numberwithin{equation}{section}

\begin{document}
\title[Divergence of the Rogers-Ramanujan Continued Fraction]
       {On the Divergence of the Rogers-Ramanujan Continued Fraction on
the Unit Circle }
\author{D. Bowman}
\address{Mathematics Department\\
       University of Illinois \\
        Champaign-Urbana, Illinois 61820}
\email{bowman@math.uiuc.edu}
\author{J. Mc Laughlin}
\address{Mathematics Department\\
       University of Illinois \\
        Champaign-Urbana, Illinois 61820}
\email{jgmclaug@math.uiuc.edu}
\thanks{The second author's research supported in part by a 
Trjitzinsky Fellowship.}
\date{April, 18, 2001}
\begin{abstract}
Let the continued fraction expansion of any irrational number $t \in (0,1)$
be denoted by $[0,a_{1}(t),a_{2}(t),\cdots]$ and
 let the $i$-th 
convergent of this continued fraction expansion be denoted by 
$c_{i}(t)/d_{i}(t)$. Let
\[ 
S=\{t \in (0,1): a_{i+1}(t) \geq 
 \phi^{d_{i}(t)}
\text{ infinitely often}\},
\]
 where $\phi = (\sqrt{5}+1)/2$. 
Let $Y_{S} =\{\exp(2 \pi i t): t \in S \}$.
It is shown that if $y \in Y_{S}$ then the 
Rogers-Ramanujan continued fraction, $R(y)$, diverges at $y$. $S$ is an
uncountable set of measure zero. It is also shown that there is 
an uncountable set of points, $G \subset Y_{S}$, 
such that if $y \in G$, then $R(y)$ does not converge generally.

It is further shown that $R(y)$ does not converge generally for 
$|y| > 1$. However we show that $R(y)$ does converge generally if $y$ is a 
primitive $5m$-th root of unity, some $m \in \mathbb{N}$ so that 
using a theorem of I. Schur, it converges generally at all roots of unity.  
\end{abstract}

\maketitle
\section{Introduction }

The Rogers-Ramanujan continued fraction, $R(x)$, is defined 
 as follows: 
\begin{equation*}
 R(x):= \cfrac{x^{\frac{1}{5}}}
{1+\cfrac{x}{1 + \cfrac{x^{2}}{1 + \cfrac{x^{3}}{1+
\ddots 
}}}}. 
\end{equation*}

Put $K(x)=x^{\frac{1}{5}}/R(x)$.  
This continued fraction seems to have been 
first investigated by L.J. Rogers in 1894 (\cite{R94}) and rediscovered 
by Ramanujan, sometime before 1913. 

It is an easy consequence
of Worpitsky's theorem (see \cite{LW92}) that  
$R(x)$ converges to values in $\hat{\mathbb{C}}$
for any $x$ inside the 
unit circle. In fact, many explicit evaluations of $R(e^{-\pi \sqrt{n}})$ 
and $R(-e^{-\pi \sqrt{n}})$ have been given for 
$n \in \mathbb{Q}^{+}$ (see, for example, \cite{BCC96},
  \cite{BC95}, \cite{R84} and \cite{Y01}), some of which were
asserted by Ramanujan without proof.

 It was stated by Ramanujan, and proved in \cite{ABJL92}, that
if $|x|>1$ then the odd and even convergents tend to different 
limits. 

This leaves the question of convergence on the unit circle.
Schur showed in \cite{S17} that if $x$ is a primitive $m$-th
root of unity, where $m \equiv 0$ $(\text{mod}\,\,5)$, then $K(x)$ diverges and if 
 $x$ is a primitive $m$-th
root of unity, $m \not \equiv 0 (\text{mod}\,\,5)$, then $K(x)$ converges 
and 
\begin{equation}\label{E:ScM}
K(x) = 
\lambda x^{\frac{1 - \lambda \sigma  m }{5}} K(\lambda),
\end{equation}
where $\lambda = \left( \frac{m}{5} \right)$ (the Legendre symbol)
and $\sigma$ is the least positive residue of $m \,(\text{mod}\,\, 5)$. 
Note that 
$K(1)= \phi = (\sqrt{5}+1)/2$, and $K(-1) = 1/ \phi$.
It follows that $R(x)$ takes only ten possible values at roots 
of unity. For later use we define 
\begin{equation}\label{E:rj}
R_{j} = 
\begin{cases}-\phi  \exp \displaystyle{ (2 \pi i j/5)},
&\,\, 1 \leq j \leq 5,\\
	\frac{\displaystyle{\exp(2 \pi i j/5)}}
{\displaystyle{\phi}},&\,\,6 \leq j \leq 10.
\end{cases}
\end{equation}

\vspace{5pt}

Remark: $\{R_{j}\}_{j=1}^{10}$ consists of the ten values 
taken by $R(x)$ at roots of unity.

Since Schur's paper it has been an open problem whether $K(x)$ converges 
or diverges at any point $x$ on the unit circle which is not a 
primitive root of unity.

Let the regular continued fraction expansion of any irrational $t \in (0,1)$
be denoted by $[0,a_{1}(t),a_{2}(t),\cdots]$.
 Let the $i$-th 
convergent of this continued fraction expansion be denoted by 
$c_{i}(t)/d_{i}(t)$. Occasionally we write $a_{i}$ for $a_{i}(t)$, 
$d_{i}$ for $d_{i}(t)$ etc, for simplicity, if there is no danger of 
ambiguity.
In this paper we prove the following theorem.
\begin{theorem}\label{T:t1}
Let
\begin{equation}\label{E:t1}
S=\{t \in (0,1): a_{i+1}(t) \geq 
 \phi^{d_{i}(t)}
\text{ infinitely often}\}.
\end{equation}
 Then $S$ is an uncountable set of measure zero
and, if $t \in S$ and $y = \exp (2 \pi i t)$, then $K(y)$ diverges.
\end{theorem}
As an example of a point in $S$, we give the following corollary to
Theorem \ref{T:t1}.
\begin{corollary}\label{C:ex}
 Let $t$ be the number with continued fraction expansion
equal $[0,a_{1},a_{2}, \cdots]$, where $a_{i}$ is the integer  consisting of 
a tower of $i$ twos with an $i$ an top.
\begin{multline*}
t=[0,2,2^{\displaystyle{2^{2}}},2^{
\displaystyle{2^{\displaystyle{2^{3}}}}},\cdots]=\\
0.484848484848484848484848484848484848484848484848484848484\\
8484848484848484848
4849277885083112437522992318812011 \cdots
\end{multline*}
If $y = \exp (2 \pi i t)$ then $K(y)$ diverges.
\end{corollary}

 In \cite{J86}, Jacobson introduced the concept of general convergence 
for continued fractions. 
General convergence is defined in  \cite{LW92} as follows.

 Let the $n$-th convergent of the continued fraction 
\[
M=b_{0}+\cfrac{a_{1}}{b_{1} + \cfrac{a_{2}}{b_{2} + \cfrac{a_{3}}{b_{3}+
\ddots }}}
\]
be denoted by $A_{n}/B_{n}$ and let 
\[S_{n}(w)= 
b_{0}+\cfrac{a_{1}}{b_{1} + \cfrac{a_{2}}{b_{2} + \cfrac{a_{3}}
{\begin{matrix}
&b_{3}+
\begin{matrix}
&\phantom{as}\\
\ddots
\end{matrix}\\
&\phantom{asasda}
\displaystyle{
\frac{a_{n}}{b_{n}+w}}
\end{matrix}
 }}}
=\frac{A_{n}+wA_{n-1}}{B_{n}+wB_{n-1}}.
\]
Define 
\begin{equation}\label{E:d}
d(w,z)=\frac{|z-w|}
{\sqrt{1+|w|^{2}}\sqrt{1+|z|^{2}}},
\end{equation}
if $w$ and $z$ are both finite, and 
\[
d(w, \infty) = \frac{1}{\sqrt{1+|w|^{2}}}.
\]
\textbf{Definition:} $M$ is said to converge generally to 
$f \in \hat{\mathbb{C}}$ if there exist sequences 
$\{v_{n}\}$, $\{w_{n}\} \subset \hat{\mathbb{C}}$ such that 
$\liminf d(v_{n},w_{n})>0$ and 
\[
\lim_{n \to \infty}S_{n}(v_{n})=\lim_{n \to \infty}S_{n}(w_{n}) = f.
\]

If a continued fraction  converges 
generally, then it does, in a certain  sense,
 to the ``right'' value.
More precisely, for $n=0,1,2,\cdots$, let 
\[
h_{n} = -S_{n}^{-1}(\infty).
\] We use  the following theorem from
\cite{LW92}.
\begin{theorem} The continued fraction $b_{0} + K(a_{n}/b_{n})$
converges generally to $f$ if and only if $\lim S_{n}(u_{n}) =f$
for every sequence $\{u_{n}\}$ from  $\hat{\mathbb{C}}$ such that
{\allowdisplaybreaks
\begin{align*}
&\liminf_{n \to \infty} d(u_{n},-h_{n})>0 &\text{if } f \not = \infty,\\
&\liminf_{n \to \infty} d(u_{n},-A_{n}/A_{n-1})>0 &\text{if } f = \infty.
\end{align*}
}
\end{theorem}
This theorem in turn has the following important corollary (
also from \cite{LW92}).
\begin{corollary}
Let $b_{0}+ K(a_{n}/b_{n})$ converge generally to $f$ and to $g$.
Then $f = g$.
\end{corollary}
 
Classical convergence implies general 
convergence
(take $u_{n}=0$ and $v_{n}= \infty$, for all $n$), 
but not conversely. Thus general convergence is
 a natural extension of classical covergence.

 As Schur showed in \cite{S17}, $K(x)$ does not converge in the 
classical sense when $x$ is an $m$-th root of unity, where
$m \equiv 0 (\text{mod}\,\,5)$. However $K(x)$ can be shown to 
converge generally in this case. We have the following 
proposition.
\begin{proposition}\label{P:pr1}
If  $x$ is an $m$-th root of unity, where
$m \equiv 0 (\text{mod}\,\,5)$, then $K(x)$ converges generally.
\end{proposition}
Taking Proposition \ref{P:pr1} along with Schur's theorem shows that
$K(x)$ converges generally at any root of unity.

This suggests the question of general convergence at points on the unit
circle which are not roots of unity.  
We have the following theorem. 
\begin{theorem}\label{T:t2}
Let $t$ be any irrational in $(0,1)$ for which 
there exist two subsequences of convergents
$\{c_{f_{n}}/d_{f_{n}}\}$ and $\{c_{g_{n}}/d_{g_{n}}\}$
and 
 integers $r$, $u \in \{0,1,2,3,4\}$,
integers $s$, $v \in \{1,2,3,4\}$  such that
\begin{align}\label{e:fieq}
c_{f_{n}} &\equiv r(\text{mod}\,\,5), &c_{g_{n}} \equiv u(\text{mod}\,\,5),\\
d_{f_{n}} &\equiv s(\text{mod}\,\,5), &d_{g_{n}} \equiv v(\text{mod}\,\,5).
\notag
\end{align}
and
\begin{align}\label{E:rcon}
a_{h_{n}+1} > 2 \pi(d_{h_{n}}+1)^{2}\phi^{d_{h_{n}}^2+2d_{h_{n}}}, 
\end{align}
for all $n$, where $h_{n} = f_{n}$ or $g_{n}$. 

Suppose further that 
\begin{align}\label{E:rsuv}
 R(\exp(2 \pi i r/s))  = R_{a} \not = R_{b}= 
  R(\exp(2 \pi i u/v)),
\end{align}
for some $a$, $b \in \{1,2,\cdots,10\}$.

 Let $S^{\diamond}$ denote the set of all 
$t \in (0,1)$ satisfying \eqref{e:fieq}, \eqref{E:rcon} and \eqref{E:rsuv}
and set 
{\allowdisplaybreaks
\begin{align}\label{G:ucon}
G = \{ \exp(2 \pi i t): t \in S^{\diamond} \}.\\
\phantom{as}\notag
\end{align}
}
Then $G$ is an uncountable set of measure zero
 such that if $y \in G$, then $K(y)$ does not converge 
generally.
\end{theorem}
Remark: It follows from \eqref{E:rcon} that $S^{\diamond} \subset S$.
Once again it is possible to give explicit examples of points
$y$ for which $K(y)$ does not converge 
generally and  in Corollary \ref{C:c33} we show that $K(y)$ does not 
converge generally for the the point $y$ in Corollary \ref{C:ex}.

An interesting question is what forms can divergence take. In
fact there are uncountably many points $y$ on the unit circle
such that $R(y)$ has subsequences of convergents tending to all 
ten of the $R_{j}$'s defined by \eqref{E:rj}.
We  prove the following proposition.
\begin{proposition}\label{P:p2} 
There exists an uncountable
subset $ G^{*} \subset Y_{S}$ such that if
  $y \in G^{*}$ then there exist ten sequences of positive integers, 
$\big\{n_{i,j}\big\}_{i=1}^{\infty}$, $1 \leq j \leq 10$, say,
 such that 
$\lim_{i \to \infty} R_{n_{i,j}}(y) = R_{j}$.
\end{proposition} 

This proposition is not strictly necessary for the proof of 
Theorem \ref{T:t2} but we find the existence of the set $G^{*}$
to be of interest. It is possible to give explicit examples of 
such points in  $ G^{*}$. We have the following corollary to 
Proposition \ref{P:p2}. 
\begin{corollary}\label{c:rr10}
Let the sequence of integers
$\{a_{i}\}_{i=1}^{\infty}$ be defined by 
\[
\{a_{1},a_{2},\cdots\}= \{ 0,2,\overline{1,2,1,0,0,1,2,1,0,2,2,4}\},
\]
where the bar indicates that the terms under it repeat infinitely often.
Let $t$ be the number with continued fraction expansion
given by 
\[
t =
[0,g_{1}+a_{1},g_{2}+a_{2},g_{3}+a_{3}, \cdots],
\] 
where $g_{i}$ is the integer  consisting of 
a tower of $i$ sixteens with an $i$ an top and the $a_{i}$'s are as 
above.
{\allowdisplaybreaks
\begin{multline*}
t=[0, 16, 16^{\displaystyle{16^{2}}}+2, 
16^{\displaystyle{16^{\displaystyle{16^{3}}}}}+1, \cdots]=\\
\phantom{as}\\
0.06249999999999999999999999999999999999999999999999999999999\\
9999999999999999999999999999999999999999999999999999999999999\\
9999999999999999999999999999999999999999999999999999999999999\\
9999999999999999999999999999999999999999999999999999999999999\\
9999999999999999999999999999999999999999999999999999999999999\\
9999999782707631005156114932594461198007415603592189407975407\\
1725266194446089419127033011861051603999 \cdots.
\end{multline*}
}
 If $y = \exp (2 \pi i t)$ then $R(y)$ has subsequences of 
convergents tending to all ten values taken by the Rogers-Ramanujan
continued fraction at roots of unity.
\end{corollary}

We also consider the question of general convergence outside the unit
circle. It was proved in \cite{ABJL92} that if $0<|x| < 1$ then 
the odd convergents of
$1/K(1/x)$ tend to
{\allowdisplaybreaks
\begin{align}\label{E:rr1}
1 - \frac{x}{1}
\begin{matrix}
\\
\, + \,
\end{matrix}
 \frac{x^2}{1}
\begin{matrix}
\\
\, - \,
\end{matrix}
 \frac{x^3}{1}
\begin{matrix}
\\
\, + \, \cdots
\end{matrix}
:= F_{1}(x)
\end{align}
}
\begin{flushleft}
while the even convergents tend to 
\end{flushleft}
{\allowdisplaybreaks
\begin{align}\label{E:rr2}
\frac{x}{1}
\begin{matrix}
\\
\, + \,
\end{matrix}
\frac{x^4}{1}
\begin{matrix}
\\
\, + \,
\end{matrix}
\frac{x^8}{1}
\begin{matrix}
\\
\, + \,
\end{matrix}
\frac{x^{12}}{1}
\begin{matrix}
\\
\, + \,\cdots
\end{matrix}
:= F_{2}(x).\\
\phantom{as} \notag
\end{align}
}
Worpitsky's theorem gives that each continued fraction does converge 
inside the unit circle to values  in $\hat{\mathbb{C}}$.
It is not clear that $F_{1}(x) \not = F_{2}(x)$  for \emph{all}
$x$ inside the unit circle but Worpitsky's theorem again gives that
$F_{1}(x) \not = F_{2}(x)$  for $|x|<1/4$ and for such $x$ we have the 
following proposition which implies that the 
Rogers-Ramanujan continued fraction
does not converge generally at $1/x$.
\begin{proposition}\label{P:pr2}
Let $C = b_{0} + K_{n=1}^{\infty}a_{n}/b_{n}$ be such that the odd and 
even convergents tend to different limits. Further suppose that there
exist positive constants $c_{1}$, $c_{2}$ and $c_{3}$ such that, 
for $i \geq 1$,  
{\allowdisplaybreaks
\begin{align}\label{con1}
\notag \\
c_{1} \leq |b_{i}| \leq c_{2},
\end{align}
}
and
{\allowdisplaybreaks
\begin{align}\label{con2}
\left|\frac{a_{2i+1}}{a_{2i}}\right| \leq c_{3}.\\
\notag
\end{align}
}
Then $C$ does not converge 
generally.
\end{proposition}

\section{ Divergence In The Classical Sense }

Let 
\[K_{n}(x) := 
1+K_{j=1}^{n}\displaystyle{\frac{x^{j}}{1}}=\frac{P_{n}(x)}{Q_{n}(x)}
\] 
denote the $n$-th convergent of 
$K(x)$ and let 
$R_{n}(x)=x^{\frac{1}{5}}/K_{n}(x)$.
 It is elementary that 
{\allowdisplaybreaks
\begin{align} \label{E:recu}
Q_{n+1}(x) = Q_{n}(x)+x^{n+1}Q_{n-1}(x).\\
\phantom{as} \notag
\end{align}
}
It can also easily be checked that if $|x|= 1$, then for $n \geq 1$,
{\allowdisplaybreaks
\begin{align}\label{E:pqs}
\phantom{as} \notag\\
|P_{n}(x)Q_{n-1}(x)-Q_{n}(x)P_{n-1}(x)|=1. \\
\phantom{as} \notag
\end{align}
}
It follows easily from the triangle inequality that, for $n \geq 2$, 
{\allowdisplaybreaks
\begin{align} \label{E:fb}
\phantom{as} \notag\\
|Q_{n}(x)| \leq F_{n+1}.\\
\phantom{as} \notag
\end{align}
}
where $\{F_{i}\}_{i=1}^{\infty}$ denotes the Fibonacci sequence defined by 
$F_{1}=F_{2}=1$ and $F_{i+1}=F_{i}+F_{i-1}$.

Suppose $\lim_{n \to \infty}P_{n}(y)/Q_{n}(y) = L \in \mathbb{C}$ for 
some $y$ on the unit circle. Then, by \eqref{E:pqs},

{\allowdisplaybreaks
\begin{align*}
\frac{1}{|Q_{n}(y)Q_{n-1}(y)|}&=
\left|\frac{P_{n}(y)}{Q_{n}(y)}-\frac{P_{n-1}(y)}{Q_{n-1}(y)}\right|\\
&\leq \left|\frac{P_{n}(y)}{Q_{n}(y)}-L\right|+
\left|\frac{P_{n-1}(y)}{Q_{n-1}(y)}-L\right|.\\
\end{align*}
}
Thus  $\lim_{n \to \infty}|Q_{n}(y)Q_{n-1}(y)| = 
\infty$. We will exhibit an uncountable set of points, of measure zero,
 on the unit circle such that if $y$ is one of these points then
 $\lim_{n \to \infty}|Q_{n}(y)Q_{n-1}(y)|\not = \infty$ so that $K(y)$ 
does not converge.

\begin{lemma}\label{L:l1}
With the notation of Theorem\ref{T:t1},
for $t \in S$, we have  
\[\left |t -\frac{ c_{i}(t)}{d_{i}(t)} \right | <
\frac{1}{d_{i}(t)^2 \phi^{d_{i}(t)}}
\] for infinitely many $i$. 
\end{lemma}

\begin{proof}
Let $i$ be one of the infinitely many integers for which
$a_{i+1}(t) \geq  \phi^{d_{i}(t)}$ and let 
 $t_{i+1} = [a_{i+1}(t),a_{i+2}(t),\cdots]$ denote the $i$-th tail of the 
continued fraction expansion for $t$. Then 
{\allowdisplaybreaks
\begin{align*}
\left |t -\frac{ c_{i}(t)}{d_{i}(t)} \right |& =
\left |\frac{t_{i+1}c_{i}(t)+c_{i-1}(t)}{t_{i+1}d_{i}(t)+d_{i-1}(t)}-
\frac{ c_{i}(t)}{d_{i}(t)} \right |\\ &= 
\frac{1}{d_{i}(t)(t_{i+1}d_{i}(t)+d_{i-1}(t))}\\
&<\frac{1}{d_{i}(t)(a_{i+1}(t)d_{i}(t)+d_{i-1}(t))}
<\frac{1}{d_{i}(t)^2 \phi^{d_{i}(t)}}. 
\end{align*}
}
\end{proof}

\begin{lemma} \label{L:dif}
Let $x$ and $y$ be two points on the unit circle.
Then, for all integers  $n \geq 0$,
 {\allowdisplaybreaks
\begin{align}\label{E:qdif}
&|Q_{n}(x)-Q_{n}(y)| \leq n^{2}\phi^{n}|x-y|.
\end{align}
}
and
{\allowdisplaybreaks
\begin{align}\label{E:pdif}
\phantom{asdaai}
|P_{n}(x)-P_{n}(y)| \leq (n+1)^2 \phi^{n+1}|x-y|.
\end{align}
} 
\end{lemma}

\vspace{5pt}

\begin{proof}
The assertions of the lemma can easily be checked for 
$n=0,1$.

Let $\beta_{i} = |Q_{i}(x)-Q_{i}(y)|$ and 
$\delta_{i} = (i+1)  F_{i}  |x-y|$. Using \eqref{E:recu}
and \eqref{E:fb} it easily follows that 
\begin{equation} \label{E:init}
\beta_{n}  \leq 
\beta_{n-1}+\beta_{n-2} + \delta_{n-1}.
\end{equation}
Iterating this last inequality downwards gives that, for 
$r = 2, \cdots, n-1$,
\begin{equation}\label{IE:ie1}
\beta_{n} \leq F_{r}\beta_{n-r+1}+F_{r-1}\beta_{n-r}+
\sum_{i=1}^{r-1}F_{i}\delta_{n-i}.
\end{equation}
The claim is true for $r=2$ by \eqref{E:init}. Suppose it is true
for $r=2,\cdots, s$. Then  
\allowdisplaybreaks{
\begin{align*}
\beta_{n} &\leq F_{s}\beta_{n-s+1}+F_{s-1}\beta_{n-s}+
\sum_{i=1}^{s-1}F_{i}\delta_{n-i}\\
&\leq F_{s}(\beta_{n-s}+\beta_{n-s-1} + \delta_{n-s})+F_{s-1}\beta_{n-s}+
\sum_{i=1}^{s-1}F_{i}\delta_{n-i}\\
&= F_{s+1}\beta_{n-s}+F_{s}\beta_{n-s-1}+
\sum_{i=1}^{s}F_{i}\delta_{n-i}
\end{align*}
}
\begin{flushleft}
and \eqref{IE:ie1} is true by induction for $2 \leq r \leq n-1$ .
\end{flushleft}

  Recall that
$\beta_{1}=0$ and $\beta_{2}=|(1+x^{2})-(1+y^{2})| \leq 2 |x-y|$.
Now in \eqref{IE:ie1} let $r=n-1$. This gives
\allowdisplaybreaks{
\begin{align*}
\beta_{n} &\leq 2 F_{n-1}|x-y|+
\sum_{i=1}^{n-2}F_{i}\delta_{n-i}=
\sum_{i=1}^{n-1}F_{i} (n-i+1)F_{n-i}|x-y|\\
&\leq\sum_{i=1}^{n-1}\phi^{n}(n-i+1)|x-y|,
\end{align*}
}
\begin{flushleft}
using the bound 
$F_{j} \leq \phi^{j}$. This last expression simplifies to 
\end{flushleft}
\allowdisplaybreaks{
\begin{align*}
&\phi^{n}|x-y|\sum_{i=2}^{n}i < n^{2}\phi^{n}|x-y|.
\end{align*}
}
\begin{flushleft}
\eqref{E:pdif} follows similarly.
\end{flushleft}
\end{proof}
To show our set has measure zero, we use the following lemma. 
\begin{lemma}\label{L:l3a}
\text{$($\cite{RS92}$)$ }
Let $f(n)>1$ for $n=1,2,\cdots$ and suppose 
$\sum_{n=1}^{\infty}1/f(n) < \infty$. 
Then the set $S^{*} = 
\{ t \in (0,1): a_{k}(t) > f(k) \text{ infinitely often}\,\,\}$ has measure
zero. 
\end{lemma}
Suppose $x$ is a primitive $m$-th root of unity. From \cite{S17}, 
one has the following table of values
 for $P_{m-2}(x)$, $P_{m-1}(x)$,   $Q_{m-2}(x)$ and $Q_{m-1}(x)$.

\begin{table}[ht]
  \begin{center}
	\begin{tabular}{| c | c | c | c | c |}
	\hline
$m$ 	 & $P_{m-2}$ 		& $P_{m-1}$     &$Q_{m-2}$      &
	$Q_{m-1}$ \\ \hline
	 & 		        &		&		&  \\
$5 \mu$  & $0$ 			& $\displaystyle{-x^{\frac{2m}{5}}
			      -x^{\frac{-2m}{5}}}$&$\displaystyle{-x^{\frac{2m}{5}}
			      		      -x^{\frac{-2m}{5}}}$& $0$ \\
	 & 		        &		&		&  \\
$5\mu+1$ & $\displaystyle{x^{\frac{1-m}{5}}}$   & $1$           &$0$ 	& $\displaystyle{x^{\frac{-1+m}{5}}} $ \\
	 & 		        &		&		&  \\
$5\mu-1$ & $\displaystyle{x^{\frac{1+m}{5}}}$   & $1$ 		&$0$ 	& $\displaystyle{x^{\frac{-1-m}{5}}} $ \\
	 & 		        &		&		&  \\
$5\mu+2$ & $\displaystyle{-x^{\frac{1+2m}{5}}}$ & $0$ 		&$1$ 	& $\displaystyle{-x^{\frac{-1-2m}{5}}} $ \\
	 & 		        &		&		&  \\
$5\mu-2$ & $\displaystyle{-x^{\frac{1-2m}{5}}}$ & $0$ 		&$1$
	& $\displaystyle{-x^{\frac{-1+2m}{5}}} $ \\
\hline
	\end{tabular}
\phantom{asdf}\\

\vspace{10pt}

	\caption{}\label{Ta:t1}
    \end{center}
\end{table}

\vspace{25pt}

\emph{Proof of Theorem \ref{T:t1}}: 
Let $t \in S$ with convergents 
$\{c_{n}/d_{n}\}_{n=0}^{\infty}$. 
Let $y= \exp(2\pi i t)$ and let 
$x_{n}=\exp(2\pi i \, c_{n}/d_{n})$.  
By Table \ref{Ta:t1},
{\allowdisplaybreaks
\begin{align}\label{IE:ie2}
\max \{ \,\,\,
|Q_{d_{n}-1}(x_{n})|,\,\,\, |Q_{d_{n}-2}(x_{n})|\,\,\} \leq 2.
\end{align}
}
Using in turn, Lemma \ref{L:dif}, the fact that chord length is shorter
than arc length,
and Lemma \ref{L:l1}, it follows that, for 
infinitely many $n$, 
{\allowdisplaybreaks
\begin{align}\label{in:1}
|Q_{d_{n}-1}(x_{n})-Q_{d_{n}-1}(y)|
&\leq (d_{n}-1)^2\phi^{d_{n}-1}|x_{n}-y|\\
&<(d_{n}-1)^2\phi^{d_{n}-1}2 \pi \left|\frac{c_{n}}{d_{n}}-t \right|
\notag \\
&< \left(\frac{d_{n}-1}{d_{n}}\right)^{2}\frac{2 \pi}{\phi} <4. \notag
\end{align}
}
Similarly,
{\allowdisplaybreaks
\begin{align}\label{in:2}
|Q_{d_{n}-2}(x_{n})-Q_{d_{n}-2}(y)|
<\left(\frac{d_{n}-2}{d_{n}}\right)^{2}\frac{2 \pi}{\phi^2} <4.
\end{align}
}
Applying the triangle inequality to \eqref{in:1} and \eqref{in:2}
and using
\eqref{IE:ie2} gives $|Q_{d_{n}-1}(y)| < 6$ and 
$|Q_{d_{n}-2}(y)| < 6$. Finally, we have that 
\[
 |Q_{d_{n}-1}(y)Q_{d_{n}-2}(y)| < 36.
\]
Since this holds for infinitely many terms of the sequence
$\{d_{n}\}_{n=0}^{\infty}$ it follows that 
$\lim_{n \to \infty}Q_{n}(y)Q_{n-1}(y) \not =  \infty$ and thus 
$K(y)$ does not converge.

We next show that $S$ has measure zero (it is clearly an uncountable set).
Let $f(i) =   \phi^{F_{i}}$, $i=1,2,\cdots$.
Then it follows that
$\sum_{i=1}^{\infty}1/f(i) < \infty$. 
Let
\[S^{*} = \{t \in (0,1): a_{i+1}(t) \geq 
 \phi^{F_{i+1}}\
\text { infinitely often}\} 
\]
 so by Lemma\eqref{L:l3a} $S^{*}$ has measure zero.

Recall
$S = \{t \in (0,1): a_{i+1}(t) \geq 
 \phi^{d_{i}(t)}
\text { infinitely often}\}$.
 $S \subset S^{*}$, since $d_{i}(t) \geq F_{i+1}$, and thus  
$S$, being a subset of a set of measure zero, has measure zero. 

\begin{flushright}
$\Box$
\end{flushright}

\emph{Proof of Corollary \ref{C:ex}:}
Denote the $i$-th convergent of the continued fraction expansion of
$t$ by $c_{i}/d_{i}$. We will show that, for $i=1,2,\cdots$,
{\allowdisplaybreaks
\begin{align}\label{E:ineq1}
a_{i+1}\geq 2^{d_{i}} > \phi^{d_{i}}.
\end{align}
}
Then $K(y)$ will diverge by Theorem \ref{T:t1}.

\begin{align}\label{E:aineq}
2^{d_{i}} &\leq a_{i+1} \Longleftrightarrow d_{i} \leq
\underbrace{2^{^{\displaystyle{.^{\displaystyle{.}
^{\displaystyle{.}^{\displaystyle{2}^{\displaystyle{i+1}}}}}}}
}}_{i\,\, \text{twos}}, \\
&\phantom{as}\notag
\end{align}
where the notation indicates that the last integer consists
of a tower of $i$ twos with an $i+1$ on top.
It can be easily checked that the second inequality 
holds for $i=1,2$. Suppose  it holds for
for $i=1,2,\cdots, r-1$. Then 
{\allowdisplaybreaks
\begin{align*}
d_{r}=a_{r}d_{r-1}+d_{r-2} &\leq
\underbrace{2^{^{\displaystyle{.^{\displaystyle{.}
^{\displaystyle{.}^{\displaystyle{2}^{\displaystyle{r}}}}}}}
}}_{r\,\, \text{twos}}
\times
\underbrace{2^{^{\displaystyle{.^{\displaystyle{.}
^{\displaystyle{.}^{\displaystyle{2}^{\displaystyle{r}}}}}}}
}}_{(r-1)\,\, \text{twos}}
+\underbrace{2^{^{\displaystyle{.^{\displaystyle{.}
^{\displaystyle{.}^{\displaystyle{2}^{\displaystyle{r-1}}}}}}}
}}_{(r-2)\,\, \text{twos}}\\
\leq\underbrace{2^{^{\displaystyle{.^{\displaystyle{.}
^{\displaystyle{.}^{\displaystyle{2}^{\displaystyle{r+1}}}}}}}
}}_{r\,\, \text{twos}}.
\end{align*}
}
Thus the first inequality in \eqref{E:aineq} holds for all
positive integers $i$ and the result follows.
\begin{flushright}
$\Box$
\end{flushright}

We will in fact show later that $K(y)$ does not converge generally
when $y$ has the value stated in the corollary above.
Note that the 
convergents of the continued fraction expansion converge very fast to 
$t$ -- the third convergent agrees with $t$ to over $10^{19700}$ 
decimal places!

Remark: It is possible to replace the set $S$ by a similar set, say
\[S_{\kappa} = \{t \in (0,1): a_{i+1}(t) \geq 
 \kappa\,\phi^{d_{i}(t)}
\text{ infinitely often}\},
\] where $\kappa$ is a positive constant and 
Theorem \ref{T:t1} will still hold for all $t$ in $S_{\kappa}$.
However $\bigcup_{\kappa \in \mathbb{R^{+}}}S_{\kappa}$ will 
still have measure zero.

\section{Divergence in the Generalized Sense}

For ease of notation later, define $Y_{S}=\{\exp(2 \pi i t): t \in S \}$.
We first prove the general convergence of $K(y)$ when $y$ is a primitive
$5m$-th root of unity, some $m \in \mathbb{N}$.

\vspace{25pt}

\emph{Proof of Proposition \ref{P:pr1}:}
From \cite{S17}, for $0 \leq r < m$,
\begin{equation}\label{E:recurs}
P_{qm+r} = P_{r}P_{m-1}^{q},\,\,\,\,\,\, Q_{qm+r}=Q_{r}Q_{m-2}^{q}.
\end{equation}
From Table \ref{Ta:t1},
\[
P_{m-1}= -x^{\frac{2m}{5}} -x^{-\frac{2m}{5}},\,\,\,\,\,\,
Q_{m-2}= -x^{\frac{m}{5}} -x^{-\frac{m}{5}}.
\]
Let $\{u_{n}\}_{n=1}^{\infty}$ be a sequence in $\hat{\mathbb{C}}$.
It is convenient to split $n \in \mathbb{Z^{+}}$ into residue 
classes  modulo $m$. We put $n = q m + r$. From \eqref{E:recurs},
\begin{equation}\label{E:userecur}
S_{n}(u_{n})=\frac{P_{n}+u_{n}P_{n-1}}{Q_{n}+u_{n}Q_{n-1}}=
\begin{cases}
\displaystyle{\left(\frac{P_{m-1}}{Q_{m-2}}\right)^{q}} \displaystyle{
\frac{P_{r}+u_{n}P_{r-1}}{Q_{r}+u_{n}Q_{r-1}}},&
 1\leq r \leq m-1,\\
\phantom{asdasd}\\
\displaystyle{\left(\frac{P_{m-1}}{Q_{m-2}}\right)^{q}}  \displaystyle{
\frac{P_{m-1}+u_{n}P_{m-1}}{Q_{m-2}+u_{n}Q_{m-1}}}, &
 r = m.
\end{cases}
\end{equation}
Suppose that  $x^{\frac{m}{5}}$
 is in the second or third quadrants. Then  
\[
|P_{m-1}|=2\cos \left(\frac{2 \pi}{5}\right) <|Q_{m-2}|
=2 \cos \left(\frac{\pi}{5}\right).
\]
Hence
\begin{equation}\label{IE:ie3}
\left| \frac{P_{m-1}}{Q_{m-2}} \right| < 1.
\end{equation}
We now construct two sequences $\{v_{n}\}$ and $\{w_{n}\}$ 
which  satisfy the conditions for general convergence at $x$.
Let
\[
M= \max_{1\leq r \leq m} \left\{ \left|\frac{Q_{r}}{Q_{r-1}} \right|
:Q_{r-1} \not = 0 \right\}.
\]

\vspace{5pt}

Put $v_{n} = M+1$ and $w_{n}= M+2$, for $n=1,2, \cdots$.
Hence 
\[
\liminf d(v_{n},w_{n})>0,
\]
 and by \eqref{E:userecur} and \eqref{IE:ie3},
\[
\lim_{n \to \infty}\frac{P_{n}+v_{n}P_{n-1}}{Q_{n}+v_{n}Q_{n-1}}
=\lim_{n \to \infty}\frac{P_{n}+w_{n}P_{n-1}}{Q_{n}+w_{n}Q_{n-1}} =0.
\] 
Thus $K(x)$ converges generally to $0$ in this case.

Next suppose that  $x^{\frac{m}{5}}$
 is in the first or fourth quadrants so that 
\[
|P_{m-1}|=2\cos \left(\frac{ \pi}{5}\right) >|Q_{m-2}|
=2 \cos \left(\frac{2\pi}{5}\right).
\]
Then
\begin{equation}\label{IE:ie4}
\left| \frac{P_{m-1}}{Q_{m-2}} \right| > 1.
\end{equation}
In this case let 
\[
M= \max_{1\leq r \leq m} \left\{ \left|\frac{P_{r}}{P_{r-1}} \right|
:P_{r-1} \not = 0 \right\}.
\]

\vspace{5pt}

As before, let $v_{n} = M+1$ and $w_{n}= M+2$, for $n=1,2, \cdots$.
Hence 
\[
\liminf d(v_{n},w_{n})>0,
\]
 and  by \eqref{E:userecur} and \eqref{IE:ie4},
\[
\lim_{n \to \infty}\frac{P_{n}+v_{n}P_{n-1}}{Q_{n}+v_{n}Q_{n-1}}
=\lim_{n \to \infty}\frac{P_{n}+w_{n}P_{n-1}}{Q_{n}+w_{n}Q_{n-1}} =\infty.
\] 
Thus $K(x)$ converges generally to $\infty$ in second case.
\begin{flushright}
$\Box$
\end{flushright}

\vspace{25pt}

\emph{Proof of Proposition \ref{P:pr2}}
Let the $i$-th convergent of 
$C=b_{0} + K_{n=1}^{\infty}a_{n}/b_{n}$  be 
denoted $A_{i}/B_{i}$.
Suppose the odd convergents tend to $f_1$ and that 
the even convergents tend to $f_{2}$. Further suppose that 
$C$ converges generally to $f \in \hat{\mathbb{C}}$ and that  
$\{v_{n}\}$, $\{w_{n}\} \subset \hat{\mathbb{C}}$  are 
two sequences such that 
{\allowdisplaybreaks
\begin{align*}\\
\lim_{n \to \infty}\frac{A_{n}+v_{n}A_{n-1}}
     {B_{n}+v_{n}B_{n-1}}=
\lim_{n \to \infty}\frac{A_{n}+w_{n}A_{n-1}} 
     {B_{n}+w_{n}B_{n-1}}= f\\
\end{align*}
}
and
{\allowdisplaybreaks
\begin{align*}
\liminf_{n \to \infty} d(v_{n},w_{n}) > 0.
\end{align*}
}
It will be shown that these two conditions lead to a contradiction.
Suppose first that $|f|< \infty$ and, without loss of generality,
that $f \not = f_{1}$. (If $f = f_{1}$ then  $f \not = f_{2}$
and we proceed similarly).
We write 
{\allowdisplaybreaks
\begin{align*}\\
\frac{A_{n}+w_{n}A_{n-1}} 
     {B_{n}+w_{n}B_{n-1}}= f+ \gamma_{n},\phantom{asdsad}
\frac{A_{n}+v_{n}A_{n-1}}
     {B_{n}+v_{n}B_{n-1}}= f + \gamma_{n}^{'},\\
\end{align*}
}
where $ \gamma_{n} \to 0$ and $ \gamma_{n}^{'} \to 0$ as 
$n \to \infty$.
 By assumption we have, for $n \geq 0$, that
$A_{2n} = B_{2n}(f_2 + \alpha_{2n})$,  
$A_{2n+1} = B_{2n+1}(f_1 + \alpha_{2n+1})$, where
$\alpha_{i} \to 0$ as $i \to \infty$. Then
{\allowdisplaybreaks
\begin{align*}
\\
\frac{A_{2n}+w_{2n}A_{2n-1}} 
     {B_{2n}+w_{2n}B_{2n-1}} &=
\frac{B_{2n}(f_2 + \alpha_{2n}) +w_{2n}B_{2n-1}(f_1 + \alpha_{2n-1})} 
     {B_{2n}+w_{2n}B_{2n-1}} \\
&= f + \gamma_{2n}.\\
\end{align*}
}
By simple algebra we have 
{\allowdisplaybreaks
\begin{align*}
\\
w_{2n}=\frac{{B_{2\,n}}\,\left( -f + {f_2} + {{\alpha }_{2\,n}} - 
      {{\gamma }_{2\,n}} \right) }{{B_{ 2\,n-1}}\,
    \left( f - {f_1} - {{\alpha }_{ 2\,n-1}} + 
      {{\gamma }_{2\,n}} \right) }.
\end{align*}
}
Similarly,
{\allowdisplaybreaks
\begin{align*}
v_{2n}=\frac{{B_{2\,n}}\,\left( -f + {f_2} + {{\alpha }_{2\,n}} - 
      {{\gamma }_{2\,n}^{'}} \right) }{{B_{ 2\,n-1}}\,
    \left( f - {f_1} - {{\alpha }_{ 2\,n-1}} + 
      {{\gamma }_{2\,n}^{'}} \right) }.\\
\end{align*}
}
If $f \not = f_{2}$ then
{\allowdisplaybreaks
\begin{align*}\\
\lim_{n \to \infty}d(v_{2n},w_{2n}) \leq 
\lim_{n \to \infty}\frac{|v_{2n} - w_{2n}|}{|w_{2n}|} =0.\\
\end{align*}
}
Hence $f = f_{2}$,
{\allowdisplaybreaks
\begin{align*}
\\
w_{2n}=\frac{{B_{2\,n}}\,\left( {{\alpha }_{2\,n}} - 
      {{\gamma }_{2\,n}} \right) }{{B_{ 2\,n-1}}\,
    \left( f - {f_1} - {{\alpha }_{ 2\,n-1}} + 
      {{\gamma }_{2\,n}} \right) }
\end{align*}
}
and
{\allowdisplaybreaks
\begin{align*}
v_{2n}=\frac{{B_{2\,n}}\,( {{\alpha }_{2\,n}} - 
      {{\gamma }_{2\,n}^{'}} ) }{{B_{ 2\,n-1}}\,
    \left( f - {f_1} - {{\alpha }_{ 2\,n-1}} + 
      {{\gamma }_{2\,n}^{'}} \right) }.\\
\end{align*}
}
 Now we show  that
\[
\lim_{n \to \infty} \left|\frac{B_{2n}}{B_{2n-1}}\right| = \infty.
\]
For iff not, then there
is a sequence $\{n_{i}\}$ and a positive constant $M$ such that 
 $|B_{2n_{i}}/B_{2n_{i}-1}| \leq M$ for all $n_{i}$ and then 
{\allowdisplaybreaks
\begin{align*}\\
&\lim_{i \to \infty} d(v_{2n_{i}},w_{2n_{i}})
\leq \lim_{i \to \infty} \left| v_{2n_{i}}-w_{2n_{i}}
\right|\\
&\leq\lim_{i \to \infty} M\left|\frac{ {{\alpha }_{2\,n_{i}}} - 
      {{\gamma }_{2\,n_{i}}^{'}}  }{
     f - {f_1} - {{\alpha }_{ 2\,n_{i}-1}} + 
      {{\gamma }_{2\,n_{i}}^{'}}  }-
\frac{ {{\alpha }_{2\,n_{i}}} - 
      {{\gamma }_{2\,n_{i}}}  }{
     f - {f_1} - {{\alpha }_{ 2\,n_{i}-1}} + 
      {{\gamma }_{2\,n_{i}}}  }
\right|=0.\\
&\phantom{as}
\end{align*}
}
Similarly, 
{\allowdisplaybreaks
\begin{align*}\\
w_{2n+1} = \frac{B_{2n+1}}{B_{2n}}
\left( \frac{f_{1}-f_{2} + \alpha_{2n+1} -\gamma_{2n+1}}
{\gamma_{2n+1}-\alpha_{2n}}\right)
\end{align*}
}
and
{\allowdisplaybreaks
\begin{align*}
v_{2n+1} = \frac{B_{2n+1}}{B_{2n}}
\left( \frac{f_{1}-f_{2} + \alpha_{2n+1} -\gamma_{2n+1}^{'}}
{\gamma_{2n+1}^{'}-\alpha_{2n}}\right)\\
\end{align*}
}
We now show that
\[
\lim_{n \to \infty} \left|\frac{B_{2n+1}}{B_{2n}}\right| = 0.
\]
If not, then there
is a sequence $\{n_{i}\}$ and some $M>0$ such that 
 $|B_{2n_{i}+1}/B_{2n_{i}}|$ $\geq M$ for all $n_{i}$. Then 
$\lim_{i \to \infty} w_{2n_{i}+1}=\lim_{i \to \infty} v_{2n_{i}+1}
= \infty$ and $\lim_{i \to \infty}$ $ d(v_{2n_{i}+1},w_{2n_{i}+1})
=0$.

Finally, we show that it is impossible to have both 
$\lim_{n \to \infty}|B_{2n+1}/B_{2n}| = 0$ and 
$\lim_{n \to \infty}|B_{2n}/B_{2n-1}| = \infty$. For ease of 
notation let
$B_{n}/B_{n-1}$ be denoted by $r_{n}$, so that
$r_{2n} \to \infty$ and $r_{2n+1} \to 0$, as $n \to \infty$.  
From the recurrence relations for the $B_{i}$'s we have  
{\allowdisplaybreaks
\begin{align*}
r_{2n}(r_{2n+1}-b_{2n+1}) &= a_{2n+1}.
\end{align*}
}
and
{\allowdisplaybreaks
\begin{align*}
r_{2n-1}(r_{2n}-b_{2n}) &= a_{2n}.
\end{align*}
}
Thus
\begin{align*}
\frac{r_{2n}}{r_{2n}-b_{2n}}
&=\frac{a_{2n+1}r_{2n-1}}{a_{2n}(r_{2n+1}-b_{2n+1})},
\end{align*}
and by \eqref{con1} and \eqref{con2}
the left side tends to 1 and the right side tends to 0,
as $n \to \infty$, giving the 
required contradiction.

If $f = \infty$ then we write 
{\allowdisplaybreaks
\begin{align*}
\frac{A_{n}+w_{n}A_{n-1}} 
     {B_{n}+w_{n}B_{n-1}}=\frac{1}{\gamma_{n}},
\end{align*}
}
where $\lim_{n \to \infty}\gamma_{n} = 0$. With the $\alpha_{i}$'s
as above we find that 
{\allowdisplaybreaks
\begin{align*}\\
w_{2n} = - \frac{{B_{2\,n}}\,
      \left( -1 + {f_2}\,{{\gamma }_{2\,n}} + 
        {{\alpha }_{2\,n}}\,{{\gamma }_{2\,n}} \right) }
      {{B_{ 2\,n-1}}\,\left( -1 + 
        {f_1}\,{{\gamma }_{2\,n}} + 
        {{\alpha }_{ 2\,n+1}}\,{{\gamma }_{2\,n}} \right) }
\end{align*}
}
and
{\allowdisplaybreaks
\begin{align*}
v_{2n} = - \frac{{B_{2\,n}}\,
      \left( -1 + {f_2}\,{{\gamma }_{2\,n}^{'}} + 
        {{\alpha }_{2\,n}}\,{{\gamma }_{2\,n}^{'}} \right) }
      {{B_{ 2\,n-1}}\,\left( -1 + 
        {f_1}\,{{\gamma }_{2\,n}^{'}} + 
        {{\alpha }_{ 2\,n+1}}\,{{\gamma }_{2\,n}^{'}} \right) }.\\
\end{align*}
}
In this case it follows easily that 
$\lim_{n \to \infty}d(w_{2n}, v_{2n}) = 0$.

\begin{flushright}
$\Box$
\end{flushright}

Before proving Theorem \ref{T:t2} and Proposition \ref{P:p2}, 
it is necessary to prove some technical 
lemmas. In what follows, $x$ is a primitive $m$-th root of unity, where 
$m \not \equiv 0 (\text{mod}\,\,5)$.
 $\Bar{\phi}  = (-\sqrt{5}+1)/2$, $K_{j}=K_{j}(x)$,
$P_{j} = P_{j}(x)$ 
and $Q_{j} = Q_{j}(x)$, for $j=0,1,2, \cdots  $. Frequent use will be made
of Binet's formula for $F_k$.
\[F_{k}= \frac{\phi^{k}-\bar{\phi}^{k}}{\sqrt{5}}.
\]
Recall also that $\lim_{k \to \infty} F_{k+1}/F_{k} = \phi$.

We also use the following facts, which can be found in
 \cite{S17} or deduced from Table \ref{Ta:t1}.
{\allowdisplaybreaks
\begin{align}\label{E:Sc1}
P_{n}&=P_{m-1}P_{n-m}+P_{m-2}Q_{n-m}; \\
&\phantom{asssd} \notag \\ 
Q_{n}&=Q_{m-1}P_{n-m}+Q_{m-2}Q_{n-m}.\notag\\
&\phantom{asssd} \notag
\end{align}
}
{\allowdisplaybreaks
\begin{align}\label{E:Sc4}
P_{qm+r}&=P_{(q-1)m+r}+P_{(q-2)m+r}\\
&\phantom{asssd} \notag \\ 
Q_{qm+r}&=Q_{(q-1)m+r}+Q_{(q-2)m+r}.\notag\\
&\phantom{asssd} \notag
\end{align}
}
For $0 \leq r <m$, there exist constants  $b_{r}$ and $b_{r}^{'}$
such that
{\allowdisplaybreaks
\begin{align}\label{E:sc5}
&\phantom{asssd} \notag \\ 
Q_{qm+r} &= b_{r}\phi^{q} + b_{r}^{'}\bar{\phi}^{q}.\\
&\phantom{asssd} \notag  
\end{align}
}
{\allowdisplaybreaks
\begin{align}\label{al:val}
Q_{2m-1} &= Q_{m-1},  &P_{2m-1} = P_{m-1}+1,\\
&\phantom{asssd}\notag \\ \
 P_{2m-2} &= P_{m-2} \text{ \phantom{sasd}  {\itshape and\/} }
&Q_{2m-2} = 1+ Q_{m-2}. \notag\\
&\phantom{asssd}\notag
\end{align}
}

\begin{lemma}\label{L:l3} For $q \geq 2$,
{\allowdisplaybreaks
\begin{align}\label{L:a1}
&\,\, \phi^{q-1} \leq |Q_{qm+m-1}| \leq \phi^{q};\\
&\phantom{asd}\notag
\end{align}
}
{\allowdisplaybreaks
\begin{align}\label{L:b2}
&\,\,\,m \equiv 1, -1 (\text{mod}\,\,5)\,\,\Longrightarrow
\phi^{q-2} \leq |Q_{qm+m-2}| \leq \phi^{q-1};\\
&\phantom{asad}\notag 
\end{align}
}
{\allowdisplaybreaks
\begin{align}\label{L:c3}
&\,\,m \equiv 2, -2 (\text{mod}\,\,5)\,\, \Longrightarrow
\phi^{q} \leq |Q_{qm+m-2}| \leq \phi^{q+1}; \\
&\phantom{asssd}\notag 
\end{align}
}
{\allowdisplaybreaks
\begin{align}\label{L:d4}
&\,\,\frac{1}{\phi^{2}} \leq \left|\frac{Q_{qm+m-1}}{Q_{qm+m-2}}\right| 
\leq \phi^{2}. \\
&\phantom{asssd}\notag 
\end{align}
}
\end{lemma}
\begin{proof}
Using \eqref{E:sc5} and \eqref{al:val} it follows that 
{\allowdisplaybreaks
\begin{align*}
&\phantom{asssd}\\
Q_{qm+m-1} = \frac{Q_{m-1}}{\sqrt{5}}(\phi^{q+1}-\bar{\phi}^{q+1})
= \phi^{q+1}\frac{Q_{m-1}}{\sqrt{5}}
\left(1- \frac{(-1)^{q+1}}{\phi^{2q+2}}\right).\\
&\phantom{asssd}
\end{align*}
}
From Table \ref{Ta:t1}, $|Q_{m-1}| = 1$ and since $q \geq 2$, it follows that
\begin{align*}
&\phantom{as}\\
\frac{\phi^{q+1}}{\sqrt{5}}\left(1-\frac{1}{\phi^{6}}\right)
&\leq |Q_{qm+m-1}|
\leq \frac{\phi^{q+1}}{\sqrt{5}}\left(1+\frac{1}{\phi^{2}}\right),\\
&\phantom{as}
\end{align*}
and \eqref{L:a1} follows easily.

Applying \eqref{E:sc5} with $r = 2m-2$ and using the values from 
  \eqref{al:val} one finds  similarly  that 
{\allowdisplaybreaks
\begin{align*}
Q_{qm+m-2} =
\begin{cases}
 \frac{1}{\sqrt{5}}(\phi^{q} - \bar{\phi}^{q}), 
\,\, &m \equiv 1,-1 (\text{mod}\,\,5),\\
&\phantom{asd}\\
\frac{1}{\sqrt{5}}(\phi^{q+2} - \bar{\phi}^{q+2}) , 
\,\, &m \equiv 2,-2 (\text{mod}\,\,5).
\end{cases}\\
\end{align*}
}
\vspace{5pt}
If $m \equiv 1, -1 (\text{mod}\,\,5)$, then 
\begin{align*}
\phantom{as}\\
\frac{\phi^{q}}{\sqrt{5}}\left(1-\frac{1}{\phi^{4}}\right)
&\leq |Q_{qm+m-2}|
\leq \frac{\phi^{q}}{\sqrt{5}}\left(1+\frac{1}{\phi^{4}}\right)\\
&\phantom{as}
\end{align*}
and \eqref{L:b2} follows.
If $m \equiv 2, -2 (\text{mod}\,\,5)$, then 
\begin{align*}
&\phantom{as}\\
\frac{\phi^{q+2}}{\sqrt{5}}\left(1-\frac{1}{\phi^{8}}\right)
&\leq |Q_{qm+m-2}|
\leq \frac{\phi^{q+2}}{\sqrt{5}}\left(1+\frac{1}{\phi^{8}}\right)\\
\phantom{as}
\end{align*}
and \eqref{L:c3} follows.
 \eqref{L:d4} is an immediate consequence of the preceding inequalities. .
\end{proof}

\begin{lemma}\label{L:l4} For $q \geq 2$,
{\allowdisplaybreaks
\begin{align}\label{M:a1}
&\frac{1}{\phi^{2q+1}} \leq
\left|K_{qm+m-1}(x)-K(x)\right|
\leq \frac{1}{\phi^{2q}};
\end{align}
}
{\allowdisplaybreaks
\begin{align}\label{M:a2}
&\frac{1}{\phi^{2q-1}} \leq
\left|K_{qm+m-2}(x)-K(x)\right|
\leq \frac{1}{\phi^{2q-2}}.
\end{align}
}
{\allowdisplaybreaks
\begin{align}\label{C:c2}
\max \{|R_{qm+m-1}(x) - R(x)|,\phantom{as}
|R_{qm+m-2}(x) - R(x)| \}\leq \frac{1}{\phi^{2q-6}}.
\end{align}
}
\end{lemma}
\begin{proof}
\eqref{E:Sc4} implies that 
\[
 P_{qm+r} = F_{q}P_{m+r}+F_{q-1}P_{r}
\]
 and 
\[
Q_{qm+r} = F_{q}Q_{m+r}+F_{q-1}Q_{r}.
\]
Using \eqref{al:val} it follows that
\[
K_{qm+m-1}=
 \frac{P_{qm+m-1}}{Q_{qm+m-1}}= 
\frac{F_{q+1}P_{m-1}+F_{q}}
{F_{q+1}Q_{m-1}}.
\]
Let $ q \to \infty$ to get 
\[
K(x) = \frac{P_{m-1}\phi+1}{Q_{m-1}\phi}.
\]
Since $|Q_{m-1}|=1$ we have that 
{\allowdisplaybreaks
\begin{align*}
&\phantom{as}\\
\left|K_{qm+m-1}-K(x)\right|&=
\left|\frac{F_{q}}{F_{q+1}}-\frac{1}{\phi}\right|=
\frac{\sqrt{5}}{\phi^{2q+2}\left(1-\frac{(-1)^{q+1}}{\phi^{2q+2}}\right)}.\\
&\phantom{as}
\end{align*}
}
The last equality follows from Binet's formula.
Thus for $q\geq 2$,
{\allowdisplaybreaks
\begin{align*}
&\phantom{as}\\
\frac{\sqrt{5}}{\phi^{2q+2}\left(1+\frac{1}{\phi^{6}}\right)}
&\leq|K_{qm+m-1}-K(x)|
\leq
\frac{\sqrt{5}}{\phi^{2q+2}\left(1-\frac{1}{\phi^{6}}\right)}.\\
&\phantom{as}
\end{align*}
}
\eqref{M:a1} now follows.

Similarly, 
{\allowdisplaybreaks
\begin{align*}
\frac{P_{qm+m-2}}{Q_{qm+m-2}}=
\frac{P_{m-2}F_{q+1}}{Q_{m-2}F_{q+1}+F_{q}} \Longrightarrow
K(x) =\frac{P_{m-2}\phi}{Q_{m-2}\phi+1}.\\
&\phantom{as}
\end{align*}
}
We consider the cases $m \equiv 1, -1 (\text{mod}\,\,5)$ and 
$m \equiv 2, -2 (\text{mod}\,\,5)$ separately. 
In the first case it can be seen from Table \ref{Ta:t1} that
$Q_{m-2}=0$ and $|P_{m-2}|=1$. In this case 
\begin{align*}
\left|K_{qm+m-2}-K(x)\right|&=
\left|\frac{F_{q+1}}{F_{q}}-\phi\right|=
\frac{\sqrt{5}}{\phi^{2q}\left(1-\frac{(-1)^{q}}{\phi^{2q}}\right)}.\\
&\phantom{as}
\end{align*}
\eqref{M:a2}  follows.
For the second case it can be seen from Table \ref{Ta:t1} that
$Q_{m-2}=1$ and again $|P_{m-2}|=1$. In this case 
\begin{align*}
&\phantom{as}\\
\left|K_{qm+m-2}-K(x)\right|&=
\left|\frac{F_{q+1}}{F_{q+2}}-\frac{1}{\phi}\right|=
\frac{\sqrt{5}}{\phi^{2q+4}\left(1-\frac{(-1)^{q+2}}{\phi^{2q+4}}\right)}.\\
&\phantom{as}
\end{align*}
and \eqref{M:a2} again follows.
 \eqref{C:c2} follows from \eqref{M:a1}
and \eqref{M:a2}.
\end{proof}

\begin{lemma}\label{L:l5}
Let $q \geq 2$ and let $n = qm+m-1$ or $qm+m-2$. Let $y$ be another point 
on the unit circle.
Suppose $P_{n}(y)= P_{n}(x) + \epsilon_{1}$,
$Q_{n}(y)= Q_{n}(x) + \epsilon_{2}$, with 
$\epsilon = \max \{|\epsilon_{1}|,|\epsilon_{2}|\} < 1/2$. Then
\begin{equation}\label{E:kns}
\left| K_{n}(y) - K_{n}(x) \right| \leq \frac{10 \epsilon}{\phi^{q-2}}.
\end{equation}
If $q\geq 3$ and the angle between $x$ and $y$ 
(measured from the origin) is less than $5 \pi /3$ and $\epsilon
\leq 1/(20\phi^2)$, then
\begin{equation}\label{E:rns}
|R_{n}(y)-R_{n}(x)| < 3\phi|x-y| + \frac{60 \epsilon}{ \phi^{q-4}},
\end{equation}
and
\begin{equation}\label{E:rnr}
|R_{n}(y) - R(x)| \leq 3\phi|x-y| + \frac{60 \epsilon}{ \phi^{q-4}} + 
\frac{1}{\phi^{2q-3}}.
\end{equation}
\end{lemma}

\begin{proof}
{\allowdisplaybreaks
\begin{align*}
\left| K_{n}(y) - K_{n}(x) \right| 
&=\left|\frac{P_{n}(y)}{Q_{n}(y)}-\frac{P_{n}(x)}{Q_{n}(x)}\right|
=\left|\frac{\epsilon_{1}Q_{n}(x)-\epsilon_{2}P_{n}(x)}
	   {Q_{n}(x)(Q_{n}(x) + \epsilon_{2})}\right|\\
&\phantom{asd}\\
&\leq \frac{|\epsilon_{1}-\epsilon_{2}|}{|Q_{n}(x) + \epsilon_{2}|}
	+\frac{|\epsilon_{2}||P_{n}(x)-Q_{n}(x)|}
	{|Q_{n}(x)||Q_{n}(x)+\epsilon_{2}|}\\
&\phantom{asd}\\
&= \frac{|\epsilon_{1}-\epsilon_{2}|}{|Q_{n}(x) + \epsilon_{2}|}
	+\frac{|\epsilon_{2}||K_{n}(x)-1|}
	{|Q_{n}(x)+\epsilon_{2}|}\\
&\phantom{asd}\\
&\leq \frac{2\epsilon}{||Q_{n}(x)|-\epsilon|}
+\frac{\epsilon\left||K(x)|+\displaystyle{
1/\phi^{2q-2}}
+1\right|}{||Q_{n}(x)|-\epsilon|}.\\
&\phantom{asd}
\end{align*}
}
Here we have used \eqref{M:a1}, \eqref{M:a2} and the bounds on 
$\epsilon_{1}$ and $\epsilon_{2}$. 
Since $|K(x)|\leq \phi$ and $\epsilon < 1/2$, it follows that 
{\allowdisplaybreaks
\begin{align*}
\left| K_{n}(y) - K_{n}(x) \right|
&\leq \frac{2\epsilon}{||Q_{n}(x)|-1/2|}
+\frac{3\epsilon}{||Q_{n}(x)|-1/2|}\\
&\phantom{asd}\\
&=
 \frac{5\epsilon}{||Q_{n}(x)|-1/2|}\\
&\phantom{asd}\\
&\leq \frac{10\epsilon}{\phi^{q-2}}.
\end{align*}
}
The last inequality follows from \eqref{L:a1}, \eqref{L:b2}, \eqref{L:c3}. 
Similarily,
{\allowdisplaybreaks
\begin{align*}
|R_{n}(y)-R_{n}(x)|
&=\left|\frac{y^{1/5}}{K_{n}(y)}-\frac{x^{1/5}}{K_{n}(x)}\right|\\
&\phantom{asdddsdsd}\\
&= \left|\frac{K_{n}(x)(y^{1/5}-x^{1/5})
+x^{1/5}(K_{n}(x)-K_{n}(y))}{K_{n}(x)K_{n}(y)}\right|\\
&\phantom{asdddsdsd}\\
&\leq \frac{|x-y|}{|K_{n}(y)|}+
\frac{|K_{n}(x)-K_{n}(y)|}{|K_{n}(x)||K_{n}(y)|}\\
&\phantom{asdddsdsd}\\
&\leq \frac{|x-y|}{\left||K_{n}(x)|-\displaystyle{
10 \epsilon / \phi^{q-2}}\right|}+
\frac{\displaystyle{
10 \epsilon/ \phi^{q-2}}}
{|K_{n}(x)|\left||K_{n}(x)|-\displaystyle{
10\epsilon/ \phi^{q-2}}\right|}.\\
&\phantom{asdddsdsd}
\end{align*}
}
Here we have used  \eqref{E:kns} and the fact that the
bound on the angle between $x$ and $y$ implies that 
$|y^{1/5}-x^{1/5}| \leq |x-y|$.
Using \eqref{M:a1}, \eqref{M:a2} and the bound on 
$\epsilon$ it follows that
{\allowdisplaybreaks
\begin{align*}
&\phantom{asdddsdsd}\\
|R_{n}(y)-R_{n}(x)|
&\leq \frac{|x-y|}{\left||K(x)|-\displaystyle{
1/ \phi^{2q-2}-1/(2\phi^{q})}\right|}\\
&\phantom{asdd}+
\frac{10 \epsilon }
{\phi^{q-2}\left||K(x)|-\displaystyle{
1/\phi^{2q-2}}\right|
\left||K(x)|-\displaystyle{
1/\phi^{2q} -1/2\phi^{q}}\right|}.\\
&\phantom{asdddsdsd}
\end{align*}
}
Since $|K(x)| = \phi$ or $1/\phi$ it follows that
{\allowdisplaybreaks
\begin{align*}
&\phantom{asdddsdsd}\\
|R_{n}(y)-R_{n}(x)|
&\leq \frac{|x-y|\phi}{1-\displaystyle{
1/\phi^{2q-3}-1/(2\phi^{q-1})}}\\
&\phantom{asdddsdsd}+
\frac{10 \epsilon}
{ \phi^{q-4}\left(1-\displaystyle{
1/\phi^{2q-3}}\right)
\left(1-\displaystyle{
1/\phi^{2q-1} -1/(2\phi^{q-1})}\right)}\\
&\phantom{asdddsdsd}\\
&\leq \frac{|x-y|\phi}{1-\displaystyle{
1/\phi^{3}-1/(2\phi^{2})}}\\
&\phantom{asdddsdsd}+
\frac{10 \epsilon}
{ \phi^{q-4}\left(1-\displaystyle{
1/\phi^{3}}\right)
\left(1-\displaystyle{
1/\phi^{3} -1/(2\phi^{2})}\right)}\\
&\phantom{asdddsdsd}\\
&\leq 3\phi|x-y| + \frac{60 \epsilon }{ \phi^{q-4}}.
\end{align*}
}
Finally, \eqref{E:rnr} follows from \eqref{C:c2} and \eqref{E:rns}.
\end{proof}

\begin{lemma}\label{L:l7}
There exists an uncountable set of points on the unit circle such 
that if $y$ is one of these points, then there exists two increasing
sequences of integers,
$\{n_{i}\}_{i=1}^{\infty}$ and $\{m_{i}\}_{i=1}^{\infty}$ say,
such that,  
{\allowdisplaybreaks
\begin{align*}
\lim_{i \to \infty}R_{n_{i}}(y) &=\lim_{i \to \infty}R_{n_{i}-1}(y) = R_{a},\\
\lim_{i \to \infty}R_{m_{i}}(y) &=\lim_{i \to \infty}R_{m_{i}-1}(y) = R_{b},
\end{align*}
}
for some $a$, $b \in \{1,2,\cdots, 10\}$, where $a \not = b$. 
\end{lemma}

\begin{proof}
With the notation of Theorem \ref{T:t2},
let $t \in S^{\diamond}$ and   
set $y= \exp (2 \pi i t)$. Let $c_{f_{n}}/d_{f_{n}}$ be one of the 
infinitely many convergents satisfying \eqref{e:fieq} and \eqref{E:rcon}
 and set $x_{n} = \exp (2 \pi i c_{f_{n}}/d_{f_{n}})$.
Then $R(x_{n}) = R_{a}$ and 
\begin{equation}\label{E:rxydif}
|x_{n}-y| < \frac{1}{ d_{f_{n}}^{2}(d_{f_{n}}+1)^{2}\phi^{d_{f_{n}}^2+2d_{f_{n}}}}.
\end{equation}

\vspace{5pt}

For the last inequality we have used the condition on the $a_{h_{n}+1}$'s 
in \eqref{E:rcon} in the same way that
the condition on the $a_{i+1}(t)$'s in \eqref{E:t1} was used in
Lemma \ref{L:l1} and the 
fact that chord length is shorter than arc length.
Let $k= d_{f_{n}}^{2}+d_{f_{n}}-1$ or $d_{f_{n}}^{2}+d_{f_{n}}-2$.
By \eqref{E:qdif}, \eqref{E:pdif} and \eqref{E:rxydif} it follows that 
\[
|P_{k}(x)-P_{k}(y)|\leq \frac{1}{\phi^{d_{f_{n}}}} 
\]
and
{\allowdisplaybreaks
\begin{align}\label{E:rpq1}
|Q_{k}(x)-Q_{k}(y)|&\leq \frac{1}{\phi^{d_{f_{n}}}}.\\
&\phantom{asd} \notag
\end{align}
}
By \eqref{E:rnr}, with $k$ as above, $q = m = d_{f_{n}}$
and $\epsilon = 1/\phi^{d_{f_{n}}}$, it follows that 
{\allowdisplaybreaks
\begin{align}\label{E:rrnrdif}
|R_{k}(y)-R_{a}| &=|R_{k}(y)-R(x_{n})|\\
 &\leq 
\frac{3 \phi}{ d_{f_{n}}^{2}(d_{f_{n}}+1)^{2}\phi^{d_{f_{n}}^2+2d_{f_{n}}}} +
\frac{60}{ \phi^{2d_{f_{n}}-4}} +\frac{1}{\phi^{2d_{f_{n}}-3}} \notag\\
&\phantom{asd} \notag \\
&\leq \frac{500}{\phi^{2d_{f_{n}}}}. \notag
\end{align}
}
Thus 
{\allowdisplaybreaks
\begin{align}
\lim_{n \to \infty}
R_{d_{f_{n}}^{2}+d_{f_{n}}-1}(y)=\lim_{n \to \infty}
R_{d_{f_{n}}^{2}+d_{f_{n}}-2}(y) = R_{a}.
\end{align}
}
Similarily,
{\allowdisplaybreaks
\begin{align}
\lim_{n \to \infty}
R_{d_{g_{n}}^{2}+d_{g_{n}}-1}(y) =\lim_{n \to \infty}
R_{d_{g_{n}}^{2}+d_{g_{n}}-2}(y) = R_{b}.
\end{align}
}

\vspace{5pt}

It is not difficult to show that $S^{\diamond}$ is an uncountable set
and from the remark following Theorem \ref{T:t2}, 
it follows that it has measure zero. Thus 
$G = \{ \exp(2 \pi i t): t \in S^{\diamond}\} $ is an uncountable set
of measure zero.
\end{proof}

\vspace{25pt}

\emph{Proof of Proposition \ref{P:p2}}:
The proof is similar to that of Lemma \ref{L:l7}.
Let 
{\allowdisplaybreaks
\begin{align*}
\phantom{as}\\
W=\{W_{i}\}_{i=1}^{12} = \{R_{6},R_{7},R_{8},R_{9},R_{10},R_{2},R_{3},R_{4},
R_{5},R_{1},R_{8},R_{7}\}.\\
\phantom{as}
\end{align*}
}
Note that $W$ contains all ten of the values taken by the Rogers-Ramanujan
continued fraction at roots of unity.
Consider the following continued fraction:
{\allowdisplaybreaks
\begin{align}\label{E:cf1}
\phantom{as} \notag \\
\alpha = 
[0,1,3,\overline{2,3,2,1,1,2,3,2,1,3,3,5}]:= [0,a_{1},a_{2},\cdots].\\
\phantom{as} \notag 
\end{align}
}
Modulo $5$, the convergents are 
\begin{align}\label{E:con}
\phantom{as} \notag \\
\left\{\frac{0}{1},\overline{\frac{1}{1},\frac{3}{4},
\frac{2}{4},\frac{4}{1},\frac{0}{1},
\frac{4}{2},\frac{4}{3},\frac{2}{3},\frac{0}{2},\frac{2}{2},\frac{2}{4}
,\frac{3}{4}
}\right\},\\
\phantom{as} \notag 
\end{align}
where the bar indicates that, modulo $5$, the convergents repeat in 
this order.

Let $t$ be any irrational in $(0,1)$ such that,
for $i \geq 1$, the $i$-th partial quotient,
$b_{i}$, and the $i$-th convergent, $c_{i}/d_{i}$,
in its continued fraction expansion, $[0,b_{1},b_{2},\cdots ]$, 
satisfy the following conditions.
{\allowdisplaybreaks
\begin{align}\label{E:gcon}
&(i)\,\, b_{i} \equiv a_{i} (\text{mod}\,\,5), \\
&(ii)\,\,\left|t-\frac{c_{i}}{d_{i}}\right| <
\frac{1}{2 \pi d_{i}^{2}(d_{i}+1)^{2}\phi^{d_{i}^2+2d_{i}}},\notag
\end{align} 
}
\begin{flushleft}
where the $a_{i}$'s are as in equation \eqref{E:cf1}.
\end{flushleft}

Set $y= \exp (2 \pi i t)$ and let $x_{n} = \exp (2 \pi i c_{n}/d_{n})$,
so that 
\begin{equation}\label{E:xydif}
|x_{n}-y| < \frac{1}{ d_{n}^{2}(d_{n}+1)^{2}\phi^{d_{n}^2+2d_{n}}}.
\end{equation}
Here we once again have used the fact that chord length is less than
arc length. Set $r= n (\text{mod}{12})$, for 
$n>0$.
Then it can be easily checked, using
\eqref{E:ScM} and  \eqref{E:con}, that 
\begin{equation}\label{E:rw}
R(x_{n}) = 
\begin{cases}
W_r,&\,\, r \ne 0\\
W_{12}&\,\, r = 0
\end{cases}
\end{equation}
Let $k= d_{n}^{2}+d_{n}-1$ or $d_{n}^{2}+d_{n}-2$.
By \eqref{E:qdif}, \eqref{E:pdif} and \eqref{E:xydif} it follows that 
\[
|P_{k}(x_{n})-P_{k}(y)|\leq \frac{1}{\phi^{d_{n}}} 
\]
and
{\allowdisplaybreaks
\begin{align}\label{E:pq1}
|Q_{k}(x_{n})-Q_{k}(y)|&\leq \frac{1}{\phi^{d_{n}}}.\\
&\phantom{asd} \notag
 \end{align}
}
By \eqref{E:rnr}, with $k$ as above, $q = m = d_{n}$
and $\epsilon = 1/\phi^{d_{n}}$, it follows that 
{\allowdisplaybreaks
\begin{align}\label{E:rnrdif}
|R_{k}(y)-R(x_{n})|  &\leq 
\frac{3 \phi}{ d_{n}^{2}(d_{n}+1)^{2}\phi^{d_{n}^2+2d_{n}}} +
\frac{60}{ \phi^{2d_{n}-4}} +\frac{1}{\phi^{2d_{n}-3}} \\
&\phantom{asd} \notag \\
&\leq \frac{500}{\phi^{2d_{n}}}. \notag
\end{align}
}
Next, for each $j \in \{1,2,\cdots,12\}$, define a sequence of integers
$\big\{s_{i,j}\big\}_{i=1}^{\infty}$, by setting 
$s_{i,j}= d_{12(i-1)+j}^2+d_{12(i-1)+j}$. By \eqref{E:rw}, 
$R(x_{12(i-1)+j}) = W_{j}$ and so, from \eqref{E:rnrdif},
{\allowdisplaybreaks
\begin{align*}
&|R_{(s_{i,j}-1) }(y) - W_{j}| 
\leq \displaystyle{
\frac{500}{\displaystyle{
\phi^{\displaystyle{
2d_{12(i-1)+j}}}}}};\\
&|R_{(s_{i,j}-2) }(y) - W_{j}| 
\leq \displaystyle{
\frac{500}{\displaystyle{
\phi^{\displaystyle{
2d_{12(i-1)+j}}}}}}.
\end{align*}
}
It follows that
{\allowdisplaybreaks
\begin{align*}
\lim_{i \to \infty}
R_{(s_{i,j}-1) }(y) =\lim_{i \to \infty}
R_{(s_{i,j}-2) }(y) =  W_{j}.\\
\phantom{asd}
\end{align*}
}
Both results hold for $1 \leq j \leq 12$.
Since the set $W$ contains all ten of the $R_{j}$'s the result is proved
for this particular $t$.

Let $S^{'}$ denote the set of all such $t \in (0,1)$ and  
set 
$G^{*}= \{\exp (2 \pi i t): t \in S^{'}\,\,\}$.
Clearly  $G^{*} \subset Y_{S}$ and is also uncountable.
\begin{flushright}
$\Box$
\end{flushright}

\vspace{25pt}

\emph{Proof of Theorem \ref{T:t2}}: 
Let $y$ be any point in $G$, where $G$ is as defined in the 
proof of Lemma
\ref{L:l7}, and
let $t$ be the irrational in $(0,1)$ for which 
$y=\exp (2 \pi i t)$.

Suppose $R(y)$ converges generally to $f \in \hat{\mathbb{C}}$ and that  
$\{v_{n}\}$, $\{w_{n}\}$  are 
two sequences such that 
\[
\lim_{n \to \infty}\frac{P_{n}+v_{n}P_{n-1}}
     {Q_{n}+v_{n}Q_{n-1}}=
\lim_{n \to \infty}\frac{P_{n}+w_{n}P_{n-1}} 
     {Q_{n}+w_{n}Q_{n-1}}=\frac{y^{\frac{1}{5}}}{f} :=g.
\]

Suppose first that $|g| < \infty$.  
By construction there exists two infinite strictly increasing
 sequences of positive integers 
$\{n_{i}\}_{i=1}^{\infty}$, 
$\{m_{i}\}_{i=1}^{\infty}$ $\subset \mathbb{N}$ such that 
{\allowdisplaybreaks
\begin{align*}
L_{a} :=\,\,\,\,  \frac{y^{\frac{1}{5}}}{R_{a}}\,\,=\,\,
\lim_{i \to \infty}
\frac{P_{n_{i}}(y)}{Q_{n_{i}}(y)}\,\,\,\,  &=\,\,\,\, 
\lim_{i \to \infty}\frac{P_{n_{i}-1}(y)}{Q_{n_{i}-1}(y)}\,\,\,\, 
\end{align*}
}
and
{\allowdisplaybreaks
\begin{align*}
L_{b}:=\,\,\,\,  \frac{y^{\frac{1}{5}}}{R_{b}}\,\,=\,\,
\lim_{i \to \infty}
\frac{P_{m_{i}}(y)}{Q_{m_{i}}(y)}\,\,\,\,  &=\,\,\,\, 
\lim_{i \to \infty}\frac{P_{m_{i}-1}(y)}{Q_{m_{i}-1}(y)},
\end{align*}
}
for some $a \ne b$, $a,b \in \{1,2,\cdots,10\}$. 
Also by construction each $n_{i}$ has the form
$d_{k_{i}}^{2}+d_{k_{i}}-1$, where $d_{k_{i}}$ is some denominator convergent 
in the continued fraction expansion of $t$, and 
likewise for each $m_{i}$. 
 It can be further  assumed that 
$L_{a} \ne g$, since 
$L_{a} \ne L_{b}$. For ease of notation write
{\allowdisplaybreaks
\begin{align*}
 &P_{n_{i}}(y)\, =\,P_{n_{i}},
&Q_{n_{i}}(y)\, =\,Q_{n_{i}},\,\phantom{aas}\\
  &P_{n_{i}-1}(y)\, =\,P_{n_{i}-1},
&Q_{n_{i}-1}(y)\, =\,Q_{n_{i}-1}.
\end{align*}
}
Write 
$P_{n_{i}}= Q_{n_{i}}(L_{a} + \epsilon_{n_{i}})$ and
$P_{n_{i}-1}= Q_{n_{i}-1}(L_{a} + \delta_{n_{i}})$, where 
$\epsilon_{n_{i}} \to 0$ and 
$\delta_{n_{i}} \to 0$ as $i \to \infty$. 
Thus 
\[
\frac{Q_{n_{i}}(L_{a}+ \epsilon_{n_{i}})+
w_{n_{i}}Q_{n_{i}-1}(L_{a} + \delta_{n_{i}})}
{Q_{n_{i}}+w_{n_{i}}Q_{n_{i}-1}}= g + \gamma_{n_{i}},
\]
where  $\gamma_{n_{i}} \to 0$ 
 as $i \to \infty$. This last equation implies that
\[
w_{n_{i}}+\frac{Q_{n_{i}}}{Q_{n_{i}-1}} 
=\frac{Q_{n_{i}}}{Q_{n_{i}-1}} \times
\frac{\epsilon_{n_{i}}-\delta_{n_{i}}}
	{g-L_{a}+\gamma_{n_{i}}-\delta_{n_{i}}}.
\]
Because of  \eqref{L:d4}, 
 the fact that each $n_{i}$ has the form
$d_{k_{i}}^{2}+d_{k_{i}}-1$, where $d_{k_{i}}$ is some denominator convergent 
in the continued fraction expansion of $t$ and \eqref{E:rpq1},  
 it follows that $Q_{n_{i}}/ Q_{n_{i}-1}$ is absolutely 
bounded. Therefore the right hand side of the last 
equality tends to 0 as $i \to \infty$ and thus 
{\allowdisplaybreaks
\begin{align}\label{w's}
&\phantom{as} \notag \\
&w_{n_{i}}+Q_{n_{i}}/ Q_{n_{i}-1} \to 0\text{ as }n_{i} \to \infty.\\
&\phantom{as} \notag 
\end{align}
}
Note  that $|w_{n_{i}}| < \infty$ for all $i$ sufficiently large, since 
$|Q_{n_{i}}/ Q_{n_{i}-1}|<\infty$. Similarily,
{\allowdisplaybreaks
\begin{align}\label{v's}
v_{n_{i}}+Q_{n_{i}}/ Q_{n_{i}-1} \to 0\text{ as }n_{i} \to \infty.\\
&\phantom{as} \notag 
\end{align}
}
By the \eqref{w's}, \eqref{v's} and the triangle inequality 
\[
\lim_{i \to \infty}|v_{n_{i}}-w_{n_{i}}| = 0.
\]
Thus 
\[
\liminf d(v_{n},w_{n}) = 0. 
\]
Therefore $R(y)$ does not converge generally. A similar argument
holds in the case where $g$ is infinite.

Since $G$ is 
  uncountable and of measure zero, this proves the theorem.
\begin{flushright}
$\Box$
\end{flushright}

\begin{corollary}\label{C:c33}
Let $y$ be as in Corollary \ref{C:ex}. Then 
$K(y)$ does not converge generally.
\end{corollary}
\begin{proof}
Let $t \in (0,1)$ be such that $y = \exp(2 \pi i t)$.
Recall that $t=[0,a_{1},a_{2},$ $\cdots]$, where $a_{i}$ is the integer  
consisting of 
a tower of $i$ twos with an $i$ an top.
Modulo 5, the convergents in the continued fraction expansion of $t$
are 

\begin{equation}\label{E:conb}
\left\{\overline{\frac{0}{1},\frac{1}{2},\frac{1}{3},
\frac{2}{0},\frac{3}{3},\frac{0}{3},
\frac{3}{1},\frac{3}{4},\frac{1}{0},\frac{4}{4},\frac{0}{4},\frac{4}{3}
,\frac{4}{2},
\frac{3}{0},\frac{2}{2},\frac{0}{2},\frac{2}{4},\frac{2}{1},\frac{4}{0},
\frac{1}{1}
}\right\},
\end{equation}\\
where once again the bar indicates that the convergents repeat modulo 5 
in this order. In particular, there are two fractions, $r/s$ and $u/v$ say,
such that \eqref{e:fieq} and \eqref{E:rsuv} holds. 
Thus it is sufficient to show that 
{\allowdisplaybreaks
\begin{align}\label{E:rsuvcon}
\phantom{as}\notag \\
&\,\,\left|t-\frac{c_{i}}{d_{i}}\right| <
\frac{1}{2 \pi d_{i}^{2}(d_{i}+1)^{2}\phi^{d_{i}^2+2d_{i}}},
\\ \phantom{as}\notag
\end{align}
}
for all $i \geq 3$, where $c_{i}/d_{i}$ is the $i$-th convergent
in the continued fraction expansion of $t$.
In particular \eqref{E:rxydif} will hold and likewise a similar 
inequality when $f_{n}$ is replaced by $g_{n}$, where 
$\{c_{f_{n}}/d_{f_{n}}\}$ and $\{c_{g_{n}}/d_{g_{n}} \}$,
are  the two sequences of convergents
 corresponding to 
$r/s$ and $u/v$.  This in turn will ensure that $y \in G$ so that $K(y)$ 
will not converge generally by Theorem \ref{T:t2}. We will show that, for 
$i \geq 3$, 
{\allowdisplaybreaks
\begin{align}\label{I:aineq}
a_{i+1} > 16^{d_{i}^{2}}
\end{align}
}
This will be sufficient 
to prove the result. Indeed, let $t_{i+1} = 
[a_{i+1},a_{i+2},\cdots]$ denote the $i$-th tail of the 
continued fraction expansion for $t$. Then 
{\allowdisplaybreaks
\begin{align*}
&\phantom{as}\\
a_{i+1} &\geq 16^{\displaystyle{d_{i}^{2}}}
=4^{\displaystyle{2 \cdot d_{i}^{2}}}
>4^{\displaystyle{(d_{i}+1)^{2}}}\\
&\phantom{as}\\
&=2^{\displaystyle{(d_{i}+1)^{2}}}2^{\displaystyle{(d_{i}+1)^{2}}}\\
&\phantom{as}\\
&>2 \pi (d_{i}+1)^{2} \phi^{\displaystyle{d_{i}^{2}+2d_{i}}}
\Longrightarrow\\
&\phantom{as}\\
\left |t -\frac{ c_{i}}{d_{i}} \right |& =
\left |\frac{t_{i+1}c_{i}+c_{i-1}}{t_{i+1}d_{i}+d_{i-1}}-
\frac{ c_{i}}{d_{i}} \right |\\ 
&\phantom{as}\\
&= \frac{1}{d_{i}(t_{i+1}d_{i}+d_{i-1})}\\
&\phantom{as}\\
&<\frac{1}{d_{i}(a_{i+1}d_{i}+d_{i-1})}\\
&\phantom{as}\\
&<\frac{1}{d_{i}^{2}a_{i+1}}\\
&\phantom{as}\\
&< \frac{1}{2 \pi d_{i}^{2}(d_{i}+1)^{2}\phi^{\displaystyle{d_{i}^2+2d_{i}}}}.
&\phantom{as}
\end{align*} 
} 
Thus all that remains is to prove \eqref{I:aineq}.
The proof of this inequality is similar to that of 
\eqref{E:ineq1}.
{\allowdisplaybreaks
\begin{align}\label{E:abineq}
16^{\displaystyle{d_{i}^2}} \leq a_{i+1}
 \Longleftrightarrow 4\, d_{i}^{2} \leq
\underbrace{2^{^{\displaystyle{.^{\displaystyle{.}
^{\displaystyle{.}^{\displaystyle{2}^{\displaystyle{i+1}}}}}}}
}}_{i\,\, \text{twos}},
\end{align}
}
where the notation indicates that the last integer consists
of a tower of $i$ twos with an $i+1$ on top.
It can be easily checked that the second inequality 
holds for $i=3,4$. Suppose  it holds for
for $i=3,4,\cdots, r-1$. Then 
{\allowdisplaybreaks
\begin{align*}
4\,d_{r}^{2}&=4\,(a_{r}d_{r-1}+d_{r-2})^{2} 
\leq 4\,(4\,a_{r}d_{r-1}^{2}+4\,d_{r-2}^{2})^{2}\\ 
&\phantom{as}\\
&\leq 4\left(4\times
\underbrace{2^{^{\displaystyle{.^{\displaystyle{.}
^{\displaystyle{.}^{\displaystyle{2}^{\displaystyle{r}}}}}}}
}}_{r\times 2\text{'s}}
\times  
\underbrace{2^{^{\displaystyle{.^{\displaystyle{.}
^{\displaystyle{.}^{\displaystyle{2}^{\displaystyle{r}}}}}}}
}}_{(r-1)\,\, \text{twos}}
+\,\, 4 \times \underbrace{2^{^{\displaystyle{.^{\displaystyle{.}
^{\displaystyle{.}^{\displaystyle{2}^{\displaystyle{r-1}}}}}}}
}}_{(r-2)\,\, \text{twos}} \right)^{2}\\
&\phantom{as}\\
&\leq\underbrace{2^{^{\displaystyle{.^{\displaystyle{.}
^{\displaystyle{.}^{\displaystyle{2}^{\displaystyle{r+1}}}}}}}
}}_{r\,\, \text{twos}}.
&\phantom{as}
\end{align*}
}
Thus the first inequality in \eqref{E:abineq} holds for all
positive integers $i \geq 3$ and the result follows.
\end{proof}

\emph{Proof of Corollary \ref{c:rr10}:}
Showing that the $i$-th partial quotient,
$b_{i}$, and the $i$-th convergent, $ c_{i}/d_{i}$, of the 
continued fraction expansion of $t$ satisfy the conditions 
in \eqref{E:gcon}, for $i=1,2, \cdots$, 
 will ensure that $y \in G^{*}$, where $G^{*}$ is as defined 
in Proposition \ref{P:p2}.

The $b_{i}$'s satisfy the first of these conditions by construction and
so all that remains is to prove the second. By the same reasoning as used 
in the proof of Corollary \ref{C:c33}, it
is sufficient to show that 
\[
g_{i+1}\geq 16^{\displaystyle{d_{i}^{2}}}
\]  since
$b_{i+1}\geq g_{i+1}$. The details are omitted since the proof is 
almost identical, the only real difference being that
{\allowdisplaybreaks
\begin{align*}\label{E:acineq}
16^{\displaystyle{d_{i}^2}} \leq g_{i+1}
 \Longleftrightarrow  d_{i}^{2} \leq
\underbrace{16^{^{\displaystyle{.^{\displaystyle{.}
^{\displaystyle{.}^{\displaystyle{16}^{\displaystyle{i+1}}}}}}}
}}_{i\times16\text{'s}},
\end{align*}
}

\begin{flushright}
$\Box$
\end{flushright}

\section{Concluding Remarks}
The set of points on the unit circle for which the Rogers-Ramanujan 
continued fraction has 
been shown to diverge has measure zero. This still leaves open the 
question of convergence for the remaining points. At present, the 
authors do not see how to use the methods of the paper to tackle
this question.

In a later paper we will examine the 
question of convergence of other $q$-continued fractions on the 
unit circle, such as the 
G\"{o}llnitz-Gordon continued fraction and some other 
$q$-continued fractions of Ramanujan.

\allowdisplaybreaks{

}

\end{document}